\documentclass[12pt,leqno]{amsart}
\usepackage{amssymb}
\usepackage{amsmath,amssymb,color}
\numberwithin{equation}{section}

\newcommand{\la}{\lambda}
\newcommand{\va}{\varphi}
\newcommand{\ppp}{\partial}

\newcommand{\www}{\widetilde}

\newcommand{\pppa}{\partial_t^{\alpha}}

\newcommand{\R}{\mathbb{R}}

\newcommand{\C}{\mathbb{C}} 
\newcommand{\N}{\mathbb{N}}

\newcommand{\ooo}{\overline}
\newcommand{\OOO}{\Omega}

\newcommand{\sumk}{\sum_{k=0}^{\infty}}
\newcommand{\sumj}{\sum_{j=1}^d}
\newcommand{\sumij}{\sum_{i,j=1}^d}

\newcommand{\sumN}{\sum_{k=1}^N}
\newcommand{\sumjl}{\sum_{\ell=1}^j}

\newcommand{\hhalf}{\frac{1}{2}}
\newcommand{\DDD}{\mathcal{D}}

\allowdisplaybreaks

\title
[]
{
Uniqueness for inverse problem of determining fractional orders 
for time-fractional advection-diffusion equations
}

\author{
$^{1,2,3}$ M.~Yamamoto }

\thanks{
$^1$ Graduate School of Mathematical Sciences, The University
of Tokyo, Komaba, Meguro, Tokyo 153-8914, Japan \\
$^2$ Honorary Member of Academy of Romanian Scientists, 
Ilfov, nr. 3, Bucuresti, Romania \\
$^3$ Peoples' Friendship University of Russia 
(RUDN University) 6 Miklukho-Maklaya St, Moscow, 117198, Russian Federation
e-mail: {\tt myama@ms.u-tokyo.ac.jp}
}

\date{}
\begin{document}
\maketitle

\baselineskip 18pt

\begin{abstract}
We consider initial boundary value problems of time-fractional 
advection-diffusion equations with the zero Dirichlet boundary value
$\pppa u(x,t) = -Au(x,t)$, where 
$-A = \sumij \ppp_i(a_{ij}(x)\ppp_j) + \sumj b_j(x)\ppp_j + c(x)$.
We establish 
the uniqueness for an inverse problem of determining an order $\alpha$ of 
fractional derivatives by data $u(x_0,t)$ for $0<t<T$ at one point 
$x_0$ in a spatial domain $\OOO$.
The uniqueness holds even under assumption that $\OOO$ and $A$ are unknown,
provided that the initial value does not change signs and is 
not identically zero.  
The proof is based on the eigenfunction expansions of 
finitely dimensional approximating solutions, a decay estimate and the 
asymptotic expansions of the Mittag-Leffler functions for large
time.
\\
{\bf Key words.}  
fractional advection-diffusion equation, uniqueness, fractional order
\\
{\bf AMS subject classifications.}
35R30, 35R11
\end{abstract}

\section{Introduction}

Throughout this article, we assume that the spatial dimensions $d=1, 2, 3$. 
We can similarly argue for 
higher dimensions $d\ge 4$, but we need more regularity for initial values
in (1.1) and (1.5) later described.   
Let $\OOO \subset \R^d$ be a bounded domain with smooth boundary
$\ppp\OOO$ and let
$$
0< \alpha, \beta < 1.
$$
By $\pppa$ we denote the Caputo derivative:
$$
\pppa g(t) = \frac{1}{\Gamma(1-\alpha)}\int^t_0
(t-s)^{-\alpha}\frac{d}{ds}g(s) ds
$$
for $\alpha \in (0,1)$ (e.g., Podlubny \cite{Po}).

Throughout this article, we set  
$$
(-Av)(x) = \sumij \ppp_i(a_{ij}(x)\ppp_jv(x)) 
+ \sumj b_j(x)\ppp_jv(x) + c(x)v(x), \quad 
x\in \OOO,                                    \eqno{(1.1)}
$$
where $a_{ij} = a_{ji} \in C^1(\ooo{\OOO})$, $b_j \in C^1(\OOO)$, 
$1\le i,j\le d$, $c\in C(\ooo{\OOO})$ are all real-valued.
Moreover we assume that $c(x) \le 0$ for $x \in \ooo{\OOO}$, and 
there exists a constant $\sigma = \sigma(a_{ij}) > 0$ such that 
$$
\sumij a_{ij}(x)\zeta_i\zeta_j \ge \sigma(a_{ij}) \sum_{i=1}^d \zeta_i^2 
\quad \mbox{for all $x \in \ooo{\OOO}$ and $\zeta_1, ..., \zeta_d \in \R$}.   
                                         \eqno{(1.2)}
$$

We consider an initial boundary value problem for a time-fractional
advection-diffusion equation: 
$$
\left\{ \begin{array}{rl}
& \pppa u(x,t) = -Au(x,t), \quad x\in \OOO, \, 0<t<T, \\
& u\vert_{\ppp\OOO\times (0,T)} = 0, \\
& u(x,0) = a(x), \quad x \in \OOO.
\end{array}\right.
                                       \eqno{(1.3)}
$$

For $\alpha \in (0,1)$, the first equation in (1.3) is 
called a fractional advection-diffusion equation with the first-order 
term $\sumj b_j\ppp_ju$, and is a macroscopic model for 
anomalous diffusion in hetrogeneous media.    
Fractional diffusion equations are studied related also 
to diffusions in fractals 
and we refer for example to Mainardi \cite{Mai},
Metzler, Gl\"ockle and Nonnenmacher \cite{MGN}, 
Metzler and Klafter \cite{MK}, Roman and Alemany \cite{RA}.  
\\

Let $x_0 \in \OOO$ and $0<t<T$ be arbitrarily chosen.
The main subject of this article is \\
{\bf Inverse problem of determining the order $\alpha$}.\\
{\it 
Determine $\alpha$ by data $u(x_0,t)$ for $0<t<T$ for (1.3). 
}
\\

Several properties such as asymptotic behavior as $t \to \infty$ of 
solution $u$ to (1.1) depend on the fractional order $\alpha$ 
of the derivative.  It is known that $\alpha$ is an 
essential physical parameter characterizing the anomaly of diffusion.
Thus our inverse problem is important not only from the theoretical 
viewpoint but also for modelling actual anomalous advection-diffusion of 
substances such as 
contaminants by a relevant fractional diffusion equation.  

As for inverse problems of determining orders and other parameters,
we refer to  Alimov and Ashurov \cite{AA}, Ashurov and Umarov \cite{AU},  
Cheng, Nakagawa, Yamamoto and Yamazaki \cite{CNYY},
Hatano, Nakagawa, Wang and Yamamoto \cite{HNWY},
Janno \cite{J}, Janno and Kinash \cite{JKi}, Jin and Kian \cite{JK},
Krasnoschok, Pereverzyev, Siryk and Vasylyeva \cite{KPSV},
Li, Zhang, Jia and Yamamoto \cite{LZJY},
Li and Yamamoto \cite{LiYa}, Tatar, Tinaztepe and Ulusoy \cite{TTU}, 
Tatar and Ulusoy \cite{TU},
Yamamoto \cite{Y2}, \cite{Y3}, Yu, Jiang and Qi \cite{YJQ}, for example.
See Li, Liu and Yamamoto \cite{LLY} as a survey.

In \cite{KPSV}, the authors proved that for smooth 
function $\va(t)$ on $[0,T]$, the order $\alpha$ can be
uniquely determined only by $\pppa \va$ and $\va$ near $t=0$ only,
whether or not $\va(t)$ is related to a solution to 
a fractional equation equation, and applied such uniqueness to 
other inverse problem of determining an order as well as a solution to
semilinear subdiffusion equations.   In \cite{KPSV}, the continuity of
$u(x_0,t)$ and $\pppa u(x_0,t)$ at $t=0$ is essential and so 
initial value $a$ must be smoother.

We can prove the uniqueness by specified data of solution 
of even if the coefficients of $A$ are unknown, and 
we refer to \cite{CNYY} as early work, and see \cite{JK}, \cite{Y3} as
related articles.

However, these works do not consider advection terms, that is, 
assume that $b_1 = b_2 = \cdots = b_d = 0$ in $\OOO$ to study the 
uniqueness in determining the order $\alpha$.  
For such a symmetric $A$, relying on a well-known eigenfunction expansion of 
solution $u$ to (1.3), we can establish the uniqueness for the inverse
problem directly.
However, for non-symmetric $A$, such a method does not work, and 
to the best knowledge of the author,
there are no works on the uniqueness in determining the order for 
non-symmetric $A$ given by (1.1).

On the other hand, since the advection term $\sumj b_j(x)\ppp_ju$ is of 
lower-order and does not drastically change the structure of the equation, 
for non-symmetric $A$, we can naturally expect a similar uniqueness result  
to \cite{AU}, \cite{JK}, \cite{Y2}, \cite{Y3}.
The main purpose of this article is to prove that such 
a conjecture is correct for $A$ given by (1.1).
Moreover, suggested by \cite{CNYY}, \cite{JK} and \cite{Y3}, the order
$\alpha$ can be uniquely determined independently of the operator $A$ and 
the domain $\OOO$. 

For the formulation of our result, we introduce operators and domains.
Let $L^2(\OOO)$, $H^2(\OOO)$, $H^1_0(\OOO)$ denote usual Lebesgue space and 
Sobolev spaces (e.g., Adams \cite{Ad}).  By $\Vert\cdot\Vert_{L^2(\OOO)}$ 
and $(\cdot, \cdot)_{L^2(\OOO)}$,
we denote the norm and the scalar product in $L^2(\OOO)$ respectively.

We recall the Poincar\'e inequality:  there exists a constant 
$C(\OOO) > 0$, depending on a bounded domain $\OOO$ with smooth boundary 
$\ppp\OOO$ such that 
$$
C(\OOO)\int_{\OOO} \vert w\vert^2 dx \le \int_{\OOO} \vert \nabla w\vert^2
dx\quad \mbox{for $w\in H^1_0(\OOO)$}.             \eqno{(1.4)}
$$

For a bounded domain $\OOO \subset \R^d$ with smooth boundary $\ppp\OOO$,
we define $-A$ by (1.1) and assume (1.2), and for another
bounded domain $\www{\OOO}$ with smooth boundary $\ppp\www{\OOO}$, 
we additionally set
$$
(-\www{A}v)(x) = \sumij \ppp_i(\www{a_{ij}}(x)\ppp_jv(x)) 
+ \sumj \www{b_j}(x)\ppp_jv(x) + \www{c}(x)v(x), \quad 
x\in \OOO,                                    \eqno{(1.5)}
$$
where $\www{a_{ij}} = \www{a_{ji}} \in C^1(\ooo{\OOO})$, 
$\www{b_j} \in C^1(\OOO)$,$1\le i,j\le d$, $\www{c}\in C(\ooo{\OOO})$ 
are all real-valued.
Moreover we assume that $\www{c}(x) \le 0$ for $x \in \ooo{\OOO}$, and 
there exists a constant $\sigma = \sigma(\www{a_{ij}}) > 0$ such that 
$$
\sumij \www{a_{ij}}(x)\zeta_i\zeta_j 
\ge \sigma(\www{a_{ij}}) \sum_{i=1}^d \zeta_i^2 
\quad \mbox{for all $x \in \ooo{\OOO}$ and $\zeta_1, ..., \zeta_d \in \R$}.   
                                         \eqno{(1.6)}
$$
We consider an initial boundary value problem for a time-fractional
advection-diffusion equation: 
$$
\left\{ \begin{array}{rl}
& \ppp_t^{\beta} \www{u}(x,t) = -\www{A}\www{u}(x,t), \quad x\in \www{\OOO}, 
\, 0<t<T, \\
& \www{u}\vert_{\ppp\www{\OOO}\times (0,T)} = 0, \\
& \www{u}(x,0) = \www{a}(x), \quad x \in \www\OOO.
\end{array}\right.
                                       \eqno{(1.7)}
$$
We define an operator $A$ in $L^2(\OOO)$ attaching the zero Dirichlet 
boundary condition, that is, the domain $\DDD(A)$ of $A$ is 
$H^2(\OOO) \cap H^1_0(\OOO)$, and the same for $\www{A}$:
$\DDD(\www{A}) := H^2(\OOO) \cap H^1_0(\OOO)$.
Then it is known that the spectrum $\sigma(A)$ of $A$ consists entirely
of eigenvalues with finite multiplicities (e.g,, Agmon \cite{Ag}).
With some numbering, we can set
$$
\sigma(A) = \{ \la_n\}_{n\in \N} \subset \C, \quad 
\lim_{n\to\infty} \mbox{Re}\, \la_n = \infty.
$$
We note $\la_n \not\in \R$ in general.  Similarly we set 
$\sigma(\www{A}) = \{ \www{\la_n}\}_{n\in \N}$.
\\

We further asssume
$$
\left\{\begin{array}{rl}
& \inf_{x\in \OOO}\left(
\hhalf\sumj \ppp_jb_j(x) - c(x)\right) + C(\OOO)\sigma(a_{ij}) > 0, \quad 
c(x) \le 0, \quad x\in \ooo{\OOO}, \\
& \inf_{x\in \www\OOO}\left(
\hhalf\sumj \ppp_j\www{b_j}(x) - \www{c}(x)\right) 
+ C(\www{\OOO})\sigma(\www{a_{ij}}) > 0, \quad 
\www{c}(x) \le 0, \quad x\in \ooo{\www\OOO}.
\end{array}\right.
                                    \eqno{(1.8)}
$$
Condition (1.8) yields
$$
\{ \mbox{Re}\, \la_n; \, \la_n \in \sigma(A)\} > 0, \quad
\{ \mbox{Re}\, \www{\la_n}; \, \www{\la_n} \in \sigma(\www{A})\} > 0.
                                    \eqno{(1.9)}
$$
For completeness, we verify (1.9) in Appendix.
\\

In this article, for simplicity, we do not discuss the case where there exist 
eigenvalues $\la_n$ such that Re $\la_n < 0$.

Now we are reay to state the main result in this article.
\\
{\bf Theorem.}\\
{\it
Let $0<\alpha, \beta < 1$ and
$\OOO \cap \www{\OOO} \ne \emptyset$, and $x_0 \in \OOO \cap
\www{\OOO}$, $T>0$ be arbitrarily chosen.  We assume that 
each of initial values $a \in \DDD(A)$ and $\www{a} \in \DDD(\www{A})$ does 
not change
signs in $\OOO$ and $\www{\OOO}$ respectively, and 
$$
\mbox{either $a\not\equiv 0$ in $\OOO$ or $\www{a}\not\equiv 0$ 
in $\www{\OOO}$}.                           \eqno{(1.10)}
$$
If 
$$
u(x_0,t) = \www{u}(x_0,t), \quad 0<t<T,
$$
then $\alpha = \beta$.
}
\\

By Lemma 1 in Section 2, we see that $u(x_0,\cdot) \in C(0,\infty)$ and so
$u(x_0,t)$ makes sense for $t>0$.
\\ 
\vspace{0.2cm}
{\bf Remark.}\\
From the proof in Section 3, we can see the following:
\\
{\it 
Let 
$$
\mbox{each of $a$ and $\www{a}$ does not change signs in $\OOO$ and 
$\www{\OOO}$ respectively.}               \eqno{(1.11)}
$$
If $u(x_0,t) = \www{u}(x_0,t)$ for $0<t<T$ and $\alpha \ne \beta$, then 
$u(x,t) = 0$ in $\OOO\times (0,T)$ and $\www{u}(x,t) = 0$ in 
$\www{\OOO} \times (0,T)$.
}
\\

In other words, under assumption (1.11), we can conclude that 
$u(x_0,t) = \www{u}(x_0,t)$ for $0<t<T$ implies either
$$
\alpha = \beta,
$$
or
$$
u = 0 \quad \mbox{in $\OOO\times (0,T)$} \quad \mbox{and}\quad
\www{u} = 0 \quad \mbox{in $\www{\OOO} \times (0,T)$}.
$$
\\

The proof relies on the asymptotic behavior of $u$ and $\www{u}$ as 
$t \to \infty$, which is an idea similar to Sakamoto and Yamamoto \cite{SY},
Yamamoto \cite{Y3}, but by the non-symmetry of $A$ and $\www{A}$, we cannot
make use of the eigenfunction expansions of tne solutions $u$ and $\www{u}$
themselves.  Alternatively we derive asymptotic expansions of 
eigenprojections of $u$ and $\www{u}$ in view of the completeness in 
$L^2(\OOO)$ of generalized eigenfunctions for each of $A$ and $\www{A}$.

The article is composed of three sections and one appendix.
In Section 2, we prepare fundamental qualitative properties of the 
solutions and a representation formula of eigenprojections of the solution.
In Section 3, on the basis of the results in Section 2, we complete the 
proof of the main result.
\\
\section{Preliminaries}

{\bf 2.1. Qualitative properties of the solution $u$ to (1.3).}

We recall that the spatial dimensions $d$ is either $1$, $2$ or $3$.
\\
{\bf Lemma 1.}
\\
{\it 
Let $a \in L^2(\OOO)$.  Then $u \in C(\ooo{\OOO}\times (0,\infty))$ for any 
$t>0$.
}
\\
{\bf Proof.}
\\
Similarly to Gorenflo, Luchko and Yamamoto \cite{GLY}, we can prove 
$u \in C((0,\infty); H^2(\OOO))$,  Since $d\le 3$, the Sobolev embedding 
yields $H^2(\OOO)\subset C(\ooo{\OOO})$, and so we see the lemma.
$\blacksquare$
\\
\vspace{0.1cm}
{\bf Lemma 2.}
\\
{\it
Let $a\in L^2(\OOO)$ and let $x_0\in \OOO$ be fixed.
Then $u(x_0,t)$ is analytic in $t>0$.
}
\\
{\bf Proof.}
\\
We can repeatt the proof of e.g., Theorem 2.2 in Li, Imanuvilov and 
Yamamoto \cite{LIY}, and we omit the details.
\\

Next we show decay estimates of the solution $u$.
\\
{\bf Proposition 1.}
\\
{\it
(i) There exists a constant $C>0$ such that 
$$
\Vert u(\cdot,t)\Vert_{L^2(\OOO)} \le \frac{C}{t^{\alpha}}
\Vert a\Vert_{L^2(\OOO)}, \quad t>0.
$$
\\
(ii) Let $d=1,2,3$ and let $x_0\in \OOO$.  Then there exists a constant $C>0$ 
such that 
$$
\vert u(x_0,t)\vert \le \frac{C}{t^{\alpha}}\Vert Aa\Vert_{L^2(\OOO)}, 
\quad t>0
$$
for each $a \in \DDD(A)$.
}
\\

The part (i) of the proposition is well-known for the case of 
symmetric $A$ and can be proved directly by the eigenfunction expansions
(e.g., \cite{SY}).  Moreover, Vergara and Zacher \cite{VZ}
proved the same decay estimate for symmetric $A$ with time dependent 
coefficients which does not admit eigenfunction expansions, and so 
the proof requires more technicality.
See also Chapter 5 in the book Kubica, Ryszewska and Yamamoto \cite{KRY}.
For non-symmetric $A$, to the best knowledge of the author,
there are no published works and so Proposition 1 (i) can be 
an independent interest.

In Appendix we provide a sketch of
the proof of Proposition 1.
\\

{\bf 2.2. Spectral decomposition and representation of solution.}

It is sufficient to argue for the operator $A$ in $\OOO$,
because for $\www{A}$ in 
$\www{\OOO}$, we can do similarly to obtain the same results.

We recall that $A$ is defined by (1.1) with $\DDD(A) = H^2(\OOO) \cap
H^1_0(\OOO)$ and satisfies (1.2), and that 
the spectrum $\sigma(A)$ consists only of isolated eigenvalues $\la_n 
\in \C$ with $n\in \N$ and that $\infty$ is a unique accumulation point
(e.g., Agmon \cite{Ag}).  For each $n\in \N$, we take
a circle $\gamma_n$ centered at $\la_n$ with sufficiently small radius 
such that $\gamma_n$ does not enclose $\la_m$ with any $m\ne n$.
We define
$$
P_n:= \frac{1}{2\pi\sqrt{-1}}\int_{\gamma_n} (z-A)^{-1} dz    \eqno{(2.1)}
$$
and then we have
$$
P_nP_m = 
\left\{ \begin{array}{rl}
0, \quad & n\ne m,\\
P_n, \quad & n=m                
\end{array}\right.
                                  \eqno{(2.2)}
$$
(e.g., Kato \cite{Ka}).
We call $P_n$ the eigenprojection for $\la_n$ and $\va \in 
P_nL^2(\OOO)$, $\ne 0$ a generalized eigenfunction.

We see that dim $P_nL^2(\OOO) < \infty$, and we set 
$$
i_n:= \mbox{dim}\, P_nL^2(\OOO) < \infty.    \eqno{(2.3)}
$$
Then 
$$
(A-\la_n)^{i_n}P_n=0                      \eqno{(2.4)}
$$
(e.g., \cite{Ka}).  We note that if $A$ is symmetric, then $\la_n\in \R$ 
$P_nL^2(\OOO) = \mbox{Ker}\, (\la_n-A)$ and 
$i_n = \mbox{dim}\, \mbox{Ker}\, (\la_n-A)$ for all $n\in \N$.

More we set
$$
D_n:= (\la_n-A)P_n, \quad n\in \N.            \eqno{(2.5)}
$$
By (2.1) and (2.2), we prove $P_n\va = \va$ for $\va \in P_nL^2(\OOO)$,
and
$$
P_nL^2(\OOO) \subset \DDD(A), \quad AP_nL^2(\OOO) \subset P_nL^2(\OOO),
\quad D_nP_nL^2(\OOO) \subset P_nL^2(\OOO),
\quad D_n^{i_n} = 0.
                                               \eqno{(2.6)}
$$
Moreover since $0 \not\in \sigma(A)$, we see that $A^{-1}$ exists 
and 
$$
A^{-1}\va = \sum_{k=0}^{i_n-1} \frac{D_n^k}{\la_n^{k+1}}\va
\in P_nL^2(\OOO) \quad \mbox{for all $\va \in P_nL^2(\OOO)$}.                           \eqno{(2.7)}
$$
Indeed (2.5) yields 
\begin{align*}
& A\left( \sum_{k=0}^{i_n-1} \frac{D_n^k}{\la_n^{k+1}}\va\right)
= (\la_n-D_n) \left( \sum_{k=0}^{i_n-1} \frac{D_n^k}{\la_n^{k+1}}\va\right)\\
=&  \sum_{k=0}^{i_n-1} \frac{D_n^k}{\la_n^k}\va
- \sum_{k=0}^{i_n-1} \frac{D_n^{k+1}}{\la_n^{k+1}}\va
= \sum_{k=0}^{i_n-1} \frac{D_n^k}{\la_n^k}\va
- \sum_{k=1}^{i_n} \frac{D_n^k}{\la_n^k}\va
= \va - \frac{D_n^{i_n}}{\la_n^{i_n}}\va = \va,
\end{align*}
which proves (2.7).
\\

We state the completeness of the generalized eigenfunctions of $A$.
\\
{\bf Lemma 3.}
\\
{\it 
(i) For any $a\in L^2(\OOO)$, there exists a sequence 
$a_N$, $N\in \N$ such that 
$$
a_N \in \sumN P_kL^2(\OOO), \quad
\lim_{N\to\infty} \Vert a_N - a\Vert_{L^2(\OOO)} = 0.
$$
\\
(ii) For any $a\in \DDD(A)$, there exists a sequence 
$a_N$, $N\in \N$ such that 
$$
a_N \in \sumN P_kL^2(\OOO) \subset \DDD(A), \quad
\lim_{N\to\infty} \Vert Aa_N - Aa\Vert_{L^2(\OOO)} = 0.
$$
In particular,
$$
\lim_{N\to\infty} \Vert a_N - a\Vert_{C(\ooo{\OOO})} = 0.    \eqno{(2.8)}
$$
}
\\
{\bf Proof.}
\\
Since the linear subspace spanned by all the generalized 
eigenfunctions of $A$ is dense in $L^2(\OOO)$ (e.g., \cite{Ag}), we see 
part (i).  Next let $a \in \DDD(A)$.  Then $Aa \in L^2(\OOO)$.
By (i) we can choose a sequence $b_N$, $N\in \N$ such that 
$$
b_N \in \sumN P_kL^2(\OOO), \quad
\lim_{N\to\infty} \Vert b_N - Aa\Vert_{L^2(\OOO)} = 0.
$$
Since $0 \not\in \sigma(A)$, we see that $A^{-1}: L^2(\OOO) \longrightarrow
L^2(\OOO)$ exists and is bounded, and 
$$
\lim_{N\to\infty} \Vert A(A^{-1}b_N) - Aa\Vert_{L^2(\OOO)} = 0.
$$
By (2.7), we see that $a_N:= A^{-1}b_N
\in \sumN P_kL^2(\OOO) \subset \DDD(A)$ and  
$$
\lim_{N\to\infty} \Vert Aa_N - Aa\Vert_{L^2(\OOO)} = 0.
$$
Therefore, in view of $d\le 3$, the Sobolev embedding 
$\DDD(A) \subset H^2(\OOO) \subset C(\ooo{\OOO})$ implies 
(2.8).  Thus the proof of Lemma 3 is complete.
$\blacksquare$
\\

Henceforth, we set
$$
\left( \begin{array}{c}
k\\
j\\
\end{array}\right)
:= \frac{k!}{j!\, (k-j)!}, \quad 0!:= 1.
$$

Now we establish a representation formula of solution to (1.3)
with $a \in P_nL^2(\OOO)$.
\\
{\bf Lemma 4.}
\\
{\it
Let $n\in \N$ be arbitrarily fixed.  Then the solution $u_n$ to (1.3)
with $a \in P_nL^2(\OOO)$ is given by 
$$
u_n(x,t) = \sumk \frac{t^{\alpha k}}{\Gamma(\alpha k + 1)}
(-\la_n + D_n)^ka(x)
= \sumk \frac{t^{\alpha k}}{\Gamma(\alpha k + 1)}
\sum_{j=0}^k 
\left( \begin{array}{c}
k\\
j\\
\end{array}\right)
(-\la_n)^{k-j}D_n^ja(x),      \eqno{(2.9)}
$$
where the series is convergent in $C(\ooo{\OOO}\times [\delta,\infty))$
with any $\delta>0$.
}
\\
{\bf Proof.}
\\
We can interpret (2.9) as
$$
E_{\alpha,1}(-(AP_n)t^{\alpha})a 
= E_{\alpha,1}((-\la_n+D_n)t^{\alpha})a
= \sumk \frac{t^{\alpha k}}{\Gamma(\alpha k + 1)}(-\la_n+D_n)^ka, 
\quad t>0,
$$
and we can verify as follows.

We set $A_n:= AP_n$.  Then $A_n: P_nL^2(\OOO) \longrightarrow 
P_nL^2(\OOO)$.  Since dim $P_nL^2(\OOO) < \infty$, we see that 
$A_n$ is a bounded operator and 
$$
\Vert A_n\Vert_{B(P_nL^2(\OOO))} =: \rho < \infty,
$$
where $\Vert A_n\Vert_{B(P_nL^2(\OOO))}$ denotes the operator norm
of $A_n: P_nL^2(\OOO) \longrightarrow P_nL^2(\OOO)$. 
The Sobolev embedding yields
$$
\Vert a\Vert_{C(\ooo{\OOO})} \le C\Vert Aa\Vert_{L^2(\OOO)}
= C\Vert A_na\Vert_{L^2(\OOO)} \le C\rho \Vert a\Vert_{L^2(\OOO)}
$$
for all $a \in P_nL^2(\OOO)$.
Moreover
$$
\Vert A_n^ka\Vert_{C(\ooo{\OOO})} \le C\Vert A_n^{k+1}a\Vert_{L^2(\OOO)} 
\le C\Vert A_n\Vert_{B(P_nL^2(\OOO))}^{k+1}\Vert a\Vert_{L^2(\OOO)}
\le C\rho^{k+1}\Vert a\Vert_{L^2(\OOO)}.
$$
Since
$$
\pppa \left( \frac{t^{\alpha k}}{\Gamma(\alpha k + 1)}\right)
= \frac{t^{(k-1)\alpha}}{\Gamma(1+(k-1)\alpha)}, \quad k\in \N,
$$
we have
\begin{align*}
& \sumk \left\vert \pppa \left( \frac{t^{\alpha k}}{\Gamma(\alpha k + 1)}
\right) \right\vert \Vert (-A_n)^k a\Vert_{C(\ooo{\OOO})}
\le \sum_{k=1}^{\infty} 
\frac{t^{\alpha (k-1)}}{\Gamma(1+(k-1)\alpha)}
C\rho^{k+1}\Vert a\Vert_{L^2(\OOO)}\\
= & C\Vert a\Vert_{L^2(\OOO)}\rho^2
\sumk \frac{(\rho t^{\alpha})^k}{\Gamma(\alpha k + 1)} < \infty,
\end{align*}
noting that 
$$
E_{\alpha,1}(\rho t^{\alpha}) 
= \sumk \frac{(\rho t^{\alpha})^k}{\Gamma(\alpha k + 1)}
$$
and $E_{\alpha,1}(z) = \sumk \frac{z^k}{\Gamma(\alpha k + 1)}$ is an entire 
function in $z\in \C$.
Therefore, we can justify
\begin{align*}
&\pppa u_n(x,t) = \sumk \pppa\left(
\frac{t^{\alpha k}}{\Gamma(\alpha k + 1)}\right) (-A_n)^ka
= \sum_{k=1}^{\infty} \frac{t^{\alpha (k-1)}}{\Gamma(1+(k-1)\alpha)}
(-A_n)^ka \\
=& -A_n\sumk \frac{t^{\alpha k}}{\Gamma(\alpha k + 1)}
(-A_n)^ka = -A_nu_n(x,t)
\end{align*}
for $x\in \ooo{\OOO}$ and $t>0$.
We can directly verify that $u_n(x,0) = a(x)$ and $u_n(\cdot,t)
\in \DDD(A) \subset H^1_0(\OOO)$, so that $u_n$ given by (2.9) satisfies
(1.3) with $a\in P_nL^2(\OOO)$.  The proof of Lemma 4 is complete.
$\blacksquare$
\\

Next we calculate the right-hand side of (2.9).  To this end, we prove
\\
{\bf Lemma 5.}
\\
{\it
For $j\in \N$ and $1\le \ell \le j$, we have
$$
\left( \frac{d}{dt^{\alpha}}\right)^j
= \sum_{\ell=1}^j \theta_{j\ell}t^{\ell-j\alpha}
\left( \frac{d}{dt}\right)^{\ell},              \eqno{(2.10)}
$$
where $\theta_{j\ell} \in \R$ are defined by
$$
\theta_{11} = \frac{1}{\alpha}, \quad
\theta_{j+1,\ell} = 
\left\{ \begin{array}{rl}
     & \frac{\theta_{j1}(1-\alpha j)}{\alpha}, \quad \ell=1, \\
     & \frac{\theta_{j\ell}(\ell-\alpha j)}{\alpha} 
              + \frac{C_{j,\ell-1}}{\alpha}, \quad 2\le \ell \le j, \\
     & \frac{\theta_{jj}}{\alpha}, \quad \ell = j + 1
\end{array}\right.
                                     \eqno{(2.11)}
$$
for $j\ge 1$.
}
\\
Here in the case of $j=1$ in (2.11), we neglect the possibility 
$2\le \ell \le j$, that is,
$$
\theta_{2,\ell} = 
\left\{ \begin{array}{rl}
     & \frac{\theta_{11}(1-\alpha)}{\alpha}, \quad \ell=1, \\
     & \frac{\theta_{11}}{\alpha}, \quad \ell = 2.
\end{array}\right.
$$
\\
{\bf Proof.}
\\
For $j=1$, we have $\frac{d}{dt^{\alpha}}
= \frac{dt}{dt^{\alpha}}\frac{d}{dt} 
= \frac{1}{\alpha}t^{1-\alpha}\frac{d}{dt}$, and (2.10) holds with 
$\theta_{11} = \frac{1}{\alpha}$.  Let (2.10) hold for $j$.  Then
\begin{align*}
&\left(\frac{d}{dt^{\alpha}}\right)^{j+1}
= \frac{d}{dt^{\alpha}}\left(
\sum_{\ell=1}^j \theta_{j\ell} t^{\ell-j\alpha}
\left(\frac{d}{dt}\right)^{\ell}\right)
= \frac{1}{\alpha}t^{1-\alpha}
\frac{d}{dt}\left(
\sum_{\ell=1}^j \theta_{j\ell} t^{\ell-j\alpha}
\left(\frac{d}{dt}\right)^{\ell}\right)\\
= &\frac{1}{\alpha}t^{1-\alpha}
\left( \sum_{\ell=1}^j \theta_{j\ell}(\ell-\alpha j)
t^{\ell-j\alpha-1}\left(\frac{d}{dt}\right)^{\ell}
+ \sum_{\ell=1}^j \theta_{j\ell}t^{\ell-j\alpha}
\left(\frac{d}{dt}\right)^{\ell+1}\right)
\end{align*}
$$
=  \sum_{\ell=1}^j \frac{\theta_{j\ell}(\ell-\alpha j)}{\alpha}
t^{\ell-(j+1)\alpha}\left(\frac{d}{dt}\right)^{\ell}
+ \sum_{\ell=2}^{j+1} \frac{\theta_{j\ell-1}}{\alpha}
t^{\ell-(j+1)\alpha}\left(\frac{d}{dt}\right)^{\ell},
                                               \eqno{(2.12)}
$$
which proves (2.10) for $j+1$, and (2.10) is seen for $j\in \N$ by the 
induction.  The recurrence formula (2.11) follows from (2.12).
Thus the proof of Lemma 5 is complete.
$\blacksquare$
\\

Now we can show a represenatio formula f $u_n(x,t)$.
\\
{\bf Proposition 2.}
\\
{\it
$$
u_n(\cdot,t) = \sum_{j=0}^{\infty} \frac{(-1)^j}{\la_n^j\, j!}
\left( \sumjl \theta_{j\ell} E_{\alpha, 1-\ell}(-\la_nt^{\alpha})\right)
D_n^jP_na \quad \mbox{in $C(\ooo{\OOO})$ for $t \ge 0$.}
                                                           \eqno{(2.13)}
$$
}
\\
{\bf Proof.}
\\
By (2.9) and $D_n^{i_n} = 0$, the absolute convergence of the series in 
$j$ and $k$ allows us to exchange the order of the summation, and we have
$$
u_n(x,t) = \sumk \frac{t^{k\alpha}}{\Gamma(\alpha k + 1)}
\sum_{j=0}^k \left( \begin{array}{c}
k\\
j\\
\end{array}\right) (-\la_n)^{k-j}D_n^ja
                                                           \eqno{(2.14)}
$$
\begin{align*}
= & \sum_{j=0}^{\infty} \frac{D_n^ja}{j!}
\sum_{k=j}^{\infty} \frac{t^{k\alpha}}{\Gamma(\alpha k + 1)}
k(k-1)\,\cdots\, (k-j+1)(-\la_n)^{k-j}\\
= & \sum_{j=0}^{\infty} \frac{D_n^ja}{j!}t^{\alpha j}
\sum_{k=j}^{\infty} \frac{k(k-1)\,\cdots\, (k-j+1)}{\Gamma(\alpha k +1)}
(-\la_nt^{\alpha})^{k-j}
= \sum_{j=0}^{\infty} 
\frac{D_n^ja}{j!}t^{\alpha j} (-1)^j \frac{d^j}{dz^j}E_{\alpha,1}(-z)\vert
_{z=\la_nt^{\alpha}}.
\end{align*}

Now we calculate $\frac{d^j}{dz^j}E_{\alpha,1}(-z)$.
Since the series $E_{\alpha,1}(z) = \sumk \frac{z^k}{\Gamma(\alpha k +1)}$ 
is convergent uniformly in any compact set of $z\in \C$, the termwise
differentiation yields
$$
\frac{d^{\ell}}{dt^{\ell}} E_{\alpha,1}(-t^{\alpha})
= t^{-\ell}E_{\alpha,1-\ell}(-t^{\alpha}), \quad \ell\in \N
                                                     \eqno{(2.15)}
$$
by 
$$
\Gamma(\alpha k + 1) = \alpha k(\alpha k-1) \cdots (\alpha k - \ell + 1)
\Gamma(\alpha k - \ell + 1).
$$
Here we note that since $\vert \Gamma(1-\ell)\vert = \infty$ for
$\ell \in \N$, we can have
$$
E_{\alpha,1-\ell}(-t^{\alpha}) = \sum_{k=0}^{\infty} \frac{(-t^{\alpha})^k}
{\Gamma(\alpha k + 1 - \ell)} 
= \sum_{k=1}^{\infty} \frac{(-t^{\alpha})^k} {\Gamma(\alpha k + 1 - \ell)}.
$$

Consequently, in view of (2.10), we have 
\begin{align*}
& \left( \frac{d}{dt^{\alpha}}\right)^j E_{\alpha,1}(-t^{\alpha})
= \sumjl \theta_{j\ell}t^{\ell-\alpha j} \frac{d^{\ell}}{dt^{\ell}}
E_{\alpha,1}(-t^{\alpha})\\
=& \sumjl \theta_{j\ell}t^{\ell-\alpha j} t^{-\ell}
E_{\alpha,1-\ell}(-t^{\alpha})
= t^{-\alpha j}\sumjl \theta_{j\ell}E_{\alpha,1-\ell}(-t^{\alpha}).
\end{align*}
Setting $t^{\alpha} = z$, we obtain
$$
\left( \frac{d}{dz}\right)^j E_{\alpha,1}(-z)
= z^{-j}\sumjl \theta_{j\ell}E_{\alpha,1-\ell}(-z)
$$
if Re $z>0$.
Therefore,
$$
\left( \frac{d}{dz}\right)^j E_{\alpha,1}(-z)\vert_{z=\la_nt^{\alpha}}
= (\la_nt^{\alpha})^{-j}\sumjl \theta_{j\ell}
E_{\alpha,1-\ell}(-\la_nt^{\alpha}).
$$
Substituting this into (2.14), we obtain (2.13).
The proof of Proposition 2 is complete.
$\blacksquare$
\\

Finally we prove a key proposition for the proof of Theorem.
Henceforth, we fix $t_0>0$ arbitrarily, and we omit the dependency 
of the constants on $t_0>0$ and $a$.
\\
{\bf Proposition 3.}
\\
{\it
There exists a constant $C(n)>0$ such that 
$$
u_n(x,t) = \frac{A^{-1}(P_na)(x)}{\Gamma(1-\alpha)t^{\alpha}}
+ R_n(x,t), \quad t\ge t_0, \, x\in \ooo{\OOO},
$$
where 
$$
\Vert R_n(\cdot,t)\Vert_{C(\ooo{\OOO})} \le \frac{C(n)}{t^{2\alpha}},
\quad t\ge t_0.
$$
}
\\
{\bf Proof.}
\\
For $\ell \in \N$, we have asymptotics
$$
E_{\alpha,1-\ell}(-\la_nt^{\alpha})
= \frac{1}{\Gamma(1-\ell-\alpha)}\frac{1}{\la_nt^{\alpha}}
+ R_{\alpha,\ell}(-\la_nt^{\alpha}), \quad t\ge t_0,      \eqno{(2.16)}
$$
where there exists a constant $C_{\alpha, \ell}>0$ such that 
$$
\vert R_{\alpha, \ell}(-\eta)\vert \le \frac{C_{\alpha,\ell}}{\eta^2},
\quad \eta > 0.                             \eqno{(2.17)}
$$
Henceforth constants $C_{\alpha,\ell}$, $C$, $C(n)>0$, etc. depend also on 
the initial value $a$ and the order 
$\alpha$, but we omit the dependency otherwise 
we need to specify it.
As for (2.16), we refer to Section 4.5 of Chapter 4 in 
Gorenflo, Kilbas, Mainardi and Rogosin \cite{GKMR},
 formula (1.8.28) (p.43) in Kilbas, Srivastava and
Trujillo \cite{KST} or Theorem 1.4 (pp.33-34) in Podlubny \cite{Po}, and 
see also Popov and Sedletskii \cite{PS} as recent research article.

We substitite (2.16) into (2.13):
$$
u_n(x,t) = \sum_{j=0}^{\infty} \frac{1}{\la_n^j}\frac{(-1)^j}{j!}
\left( \sumjl \left( \frac{\theta_{j\ell}}{\Gamma(1-\ell-\alpha)}
\frac{1}{\la_nt^{\alpha}}
+ \theta_{j\ell}R_{\alpha,\ell}(-\la_nt^{\alpha})\right)\right)
D_n^jP_na
$$
$$
=:  \sum_{j=0}^{\infty} \frac{1}{\la_n^{j+1}}
\left( \frac{(-1)^j}{j!} \left(
\sumjl \frac{\theta_{j\ell}}{\Gamma(1-\ell-\alpha)}\right)\right)
\frac{1}{t^{\alpha}}D_n^jP_na + S_n(t^{\alpha})(x).
                                                            \eqno{(2.18)}
$$
We note that by $D_n^{i_n} = 0$ the series are finite.
Here, by (2.16), we can find a constant $C_1(n)>0$ such that 
$$
\Vert S_n(t^{\alpha})(\cdot)\Vert_{C(\ooo{\OOO})}
\le \left( 
\sum_{j=0}^{\infty} \left\vert \frac{(-1)^j}{\la_n^j j!}\right\vert
\sumjl \frac{\vert \theta_{j\ell}\vert C_{\alpha,\ell}}
{\la_n^2t^{2\alpha}}\right)
\Vert D_n^jP_na\Vert_{C(\ooo{\OOO})}
\le \frac{C_1(n)}{t^{2\alpha}}                \eqno{(2.19)}
$$
for all $t \ge t_0$.

We calculate
$$
\Phi_j:= \frac{(-1)^j}{j!}
\sumjl \frac{\theta_{j\ell}}{\Gamma(1-\ell-\alpha)}, \quad 
j\in \N.
$$
By (2.11) and $\Gamma(1-\alpha) = -\alpha\Gamma(-\alpha)$, we see
$$
\Phi_1 = \frac{(-1)^1}{1!}\frac{\theta_{11}}{\Gamma(-\alpha)}
= \frac{-1}{\alpha\Gamma(-\alpha)} = \frac{1}{\Gamma(1-\alpha)}.
                                                       \eqno{(2.20)}
$$
Now we will verify $\Phi_j=\Phi_1$ for each $j\in \N$.
By (2.11) we obtain
\begin{align*}
& \Phi_{j+1} = \frac{(-1)^{j+1}}{(j+1)!}
\sum_{\ell=1}^{j+1} \frac{\theta_{j+1,\ell}}{\Gamma(1-\ell-\alpha)}\\
=& \frac{(-1)^{j+1}}{(j+1)!}
\biggl\{ \frac{\theta_{j1}(1-\alpha j)}{\alpha}\frac{1}{\Gamma(-\alpha)}
+ \sum_{\ell=2}^j \frac{1}{\Gamma(1-\ell-\alpha)}
\left( \frac{\theta_{j\ell}(\ell-\alpha j)}{\alpha} 
  + \frac{\theta_{j,\ell-1}}{\alpha}\right)\\
+ & \frac{\theta_{jj}}{\alpha}\frac{1}{\Gamma(-j-\alpha)} \biggr\}\\
=& \frac{(-1)^{j+1}}{(j+1)!}
\biggl\{ \frac{\theta_{j1}(1-\alpha j)}{\alpha\Gamma(-\alpha)}
+ \sum_{\ell=2}^j \frac{\theta_{j\ell}(\ell-\alpha j)}
{\alpha\Gamma(1-\ell-\alpha)}\\
+& \sum_{\ell=1}^{j-1} \frac{\theta_{j\ell}}{\alpha\Gamma(-\ell-\alpha)}
+ \frac{\theta_{jj}}{\alpha}\frac{1}{\Gamma(-j-\alpha)} \biggr\}.
\end{align*}
Here we used
$$
\sum_{\ell=2}^{j} \frac{\theta_{j,\ell-1}}{\alpha\Gamma(1-\ell-\alpha)}
= \sum_{\ell=1}^{j-1} \frac{\theta_{j\ell}}{\alpha\Gamma(-\ell-\alpha)}.
$$
Using $\Gamma(1-\eta) = -\eta\Gamma(-\eta)$ for $\eta \not\in \N \cup \{0\}$, 
we obtain $(-\ell-\alpha)\Gamma(-\ell-\alpha) = \Gamma(1-\ell-\alpha)$ for 
$\ell=1,2, ...., j$, and so
$$
\frac{1}{\alpha\Gamma(-\ell-\alpha)} 
= -\frac{\ell+\alpha}{\alpha\Gamma(1-\ell-\alpha)},
$$
which means
$$
\frac{\ell-\alpha j}{\alpha\Gamma(1-\ell-\alpha)}
+ \frac{1}{\alpha\Gamma(-\ell-\alpha)}
= \frac{1}{\alpha\Gamma(1-\ell-\alpha)}(\ell-\alpha j - \ell - \alpha) 
= \frac{-(j+1)}{\Gamma(1-\ell-\alpha)}.
$$
Therefore,
\begin{align*}
& \Phi_{j+1}
= \frac{(-1)^{j+1}}{(j+1)!}
\biggl\{ \left(\frac{\theta_{j1}(1-\alpha j)}{\alpha\Gamma(-\alpha)}
+ \frac{\theta_{j1}}{\alpha\Gamma(-1-\alpha)}\right)\\
+& \sum_{\ell=2}^{j-1} \left(\frac{\theta_{j\ell}(\ell-\alpha j)}
{\alpha\Gamma(1-\ell-\alpha)} 
+ \frac{\theta_{j\ell}}{\alpha\Gamma(-\ell-\alpha)}\right)
+ \left( \frac{\theta_{jj}(j-\alpha j)}{\alpha\Gamma(1-j-\alpha)}
+ \frac{\theta_{jj}}{\alpha\Gamma(-j-\alpha)}\right) \biggr\}\\
=& \frac{(-1)^{j+1}}{(j+1)!}
 \sum_{\ell=1}^j \left(
\frac{\ell-\alpha j}{\alpha\Gamma(1-\ell-\alpha)} 
+ \frac{1}{\alpha\Gamma(-\ell-\alpha)}\right)\theta_{j\ell}\\
=& \frac{(-1)^{j+1}}{(j+1)!}(-1)(j+1)
 \sum_{\ell=1}^j \frac{\theta_{j\ell}}{\Gamma(1-\ell-\alpha)}
= \frac{(-1)^j}{j!}
\sum_{\ell=1}^j \frac{\theta_{j\ell}}{\Gamma(1-\ell-\alpha)}
= \Phi_j
\end{align*} 
for each $j\in \N$.  Thus (2.20) yields 
$\Phi_j = \Phi_1 = \frac{1}{\Gamma(1-\alpha)}$ for all $j \in \N$.
\\

Hence, applying also (2.7), we obtain
\begin{align*}
& \sum_{j=0}^{\infty} \frac{1}{\la_n^{j+1}}
\left( \frac{(-1)^j}{j!}
\sumjl \frac{\theta_{j\ell}}{\Gamma(1-\ell-\alpha)}\right)
\frac{1}{t^{\alpha}}D_n^jP_na
= \sum_{j=0}^{\infty} \frac{1}{\la_n^{j+1}}\Phi_j
\frac{1}{t^{\alpha}}D_n^jP_na\\
=& \frac{1}{\Gamma(1-\alpha)} \left(\sum_{j=0}^{i_n-1} 
\frac{1}{\la_n^{j+1}}D_n^jP_na\right)\frac{1}{t^{\alpha}}
= \frac{1}{\Gamma(1-\alpha)}A^{-1}(P_na)(x)\frac{1}{t^{\alpha}}.
\end{align*}
Combining (2.19), we complete the proof of Proposition 3.
$\blacksquare$
\section{Completion of Proof of Theorem}

{\bf First Step.}

We recall that $u$ and $\www{u}$ are the solutions to 
(1.3) and (1.7) with the initial values $a$ and $\www{a}$
respectively.  For $a\in \DDD(A)$, by Lemma 3, for each $n\in \N$,
we can find $a_N \in \sum_{k=1}^N P_kL^2(\OOO)$ such that 
$$
\lim_{N\to\infty} \rho_N = 0,                       \eqno{(3.1)}
$$
where $\rho_N := \Vert A(a-a_N)\Vert_{L^2(\OOO)}$.
Let $v_N$ and $w_N$ satisfy
$$
\left\{ \begin{array}{rl}
& \pppa v_N = -Av_N \quad \mbox{in $\OOO\times (0,\infty)$}, \\
& v_N\vert_{\ppp\OOO\times (0,\infty)} = 0, \\
& v_N(x,0) = a_N(x), \quad x \in \OOO
\end{array}\right.
                                         \eqno{(3.2)}
$$
and
$$
\left\{ \begin{array}{rl}
& \pppa w_N = -Aw_N \quad \mbox{in $\OOO\times (0,\infty)$}, \\
& w_N\vert_{\ppp\OOO\times (0,\infty)} = 0, \\
& w_N(x,0) = a(x) - a_N(x), \quad x \in \OOO.
\end{array}\right.
                                         \eqno{(3.3)}
$$
Then $u = v_N + w_N$ in $\OOO \times (0,\infty)$.

First applying Proposition 1 to (3.3) and setting
$W_N(t) = w_N(x_0,t)$, we obtain
$$
\vert W_N(t)\vert \le \frac{C}{t^{\alpha}}\rho_N, \quad t\ge t_0.
                                                \eqno{(3.4)}
$$
Since $v_N = \sum_{n=1}^N u_n$ by the uniqueness of solution to (3.2), 
the application of Proposition 3 to $u_n$ yields
\begin{align*}
& v_N(x_0,t) = \sum_{n=1}^N u_n(x_0,t) 
= \frac{1}{\Gamma(1-\alpha)t^{\alpha}} \sum_{n=1}^N (A^{-1}(P_na))(x_0)
+ \sum_{n=1}^N R_n(x_0, t)\\
=& \frac{1}{\Gamma(1-\alpha)t^{\alpha}} (A^{-1}a_N)(x_0)
+ \sum_{n=1}^N R_n(x_0,t), \quad t\ge t_0.
\end{align*}
Therefore, we can find a constant $C_0 = C_0(N) > 0$ and 
a function $S_N(t)$ such that 
$$
v_N(x_0,t) = \frac{1}{\Gamma(1-\alpha)t^{\alpha}} (A^{-1}a_N)(x_0)
+ S_N(t), \quad t\ge t_0,                       \eqno{(3.5)}
$$
where 
$$
\vert S_N(t)\vert \le \frac{C_0(N)}{t^{2\alpha}}, \quad t\ge t_0.
                                                    \eqno{(3.6)}
$$
Consequently, 
$$
u(x_0,t) = \frac{1}{\Gamma(1-\alpha)t^{\alpha}} (A^{-1}a_N)(x_0)
+ S_N(t) + W_N(t), \quad t\ge t_0.                    \eqno{(3.7)}
$$

We can similarly argue for $\www{A}$ to see that there exist
$\www{a_N} \in \DDD(\www{A})$ and functions $\www{S_N}$ and 
$\www{W_N}$ for $N\in \N$ such that 
$$
\www{\rho_N}:= \Vert \www{A}(\www{a} - \www{a_N})\Vert_{L^2(\www{\OOO})}
\longrightarrow 0 \quad \mbox{as $N\to \infty$},
$$
$$
\www{u}(x_0,t) = \frac{1}{\Gamma(1-\beta)t^{\beta}} (\www{A}^{-1}\www{a_N})
(x_0) + \www{S_N}(t) + \www{W_N}(t), \quad t\ge t_0             \eqno{(3.8)}
$$
and
$$
\vert \www{S_N}(t)\vert \le \frac{C_0(N)}{t^{2\beta}}, \quad
\vert \www{W_N}(t)\vert \le \frac{C}{t^{\beta}}\www{\rho_N}, \quad t \ge t_0.
                                                                 \eqno{(3.9)}   $$

{\bf Second Step.}

By Lemma 2, from $u(x_0,t) = \www{u}(x_0,t)$ for $0<t<T$, we derive 
$$
u(x_0,t) = \www{u}(x_0,t) \quad \mbox{for all $t>0$.}
$$
Therefore, (3.7) and (3.8) yield
\begin{align*}
& \frac{1}{\Gamma(1-\alpha)t^{\alpha}} (A^{-1}a_N)(x_0) 
+ S_N(t) + W_N(t)\\
=& \frac{1}{\Gamma(1-\beta)t^{\beta}} (\www{A}^{-1}\www{a_N})
(x_0) + \www{S_N}(t) + \www{W_N}(t), \quad t\ge t_0.
\end{align*}

Assume that $\alpha < \beta$.  Then we multiply by $t^{\alpha}$, we have
$$
\frac{1}{\Gamma(1-\alpha)} (A^{-1}a_N)(x_0) 
+ t^{\alpha}S_N(t) + t^{\alpha}W_N(t)
$$
$$
= \frac{1}{\Gamma(1-\beta)t^{\beta-\alpha}} (\www{A}^{-1}\www{a_N})
(x_0) + t^{\alpha}\www{S_N}(t) + t^{\alpha}\www{W_N}(t), \quad t\ge t_0
                                         \eqno{(3.10)}
$$
for all $N\in \N$ and all $t\ge t_0$.

By (3.4), (3.6) and (3.9), for each fixed $N\in \N$, we have
$$
\left\{ \begin{array}{rl}
& \lim_{t\to\infty} \vert t^{\alpha}S_N(t)\vert 
\le \lim_{t\to\infty} \frac{C_0(N)}{t^{\alpha}} = 0, \\
& \lim_{t\to\infty} \vert t^{\alpha}\www{S_N}(t)\vert 
\le \lim_{t\to\infty} \frac{C_0(N)}{t^{2\beta-\alpha}} = 0
\end{array}\right.                                        
                              \eqno{(3.11)}
$$
by $\beta > \alpha$.  Moreover, whenever we fix $N\in \N$ arbitrarily, 
by (3.4) and (3.9) we see
$$
\vert t^{\alpha}W_N(t)\vert \le C\rho_N, \quad
\vert t^{\alpha}\www{W_N}(t)\vert \le \frac{C\www{\rho_N}}{t^{\beta-\alpha}}
\longrightarrow 0 \quad \mbox{as $t\to\infty$} 
                                                \eqno{(3.12)}
$$
and
$$
\lim_{t\to\infty} \frac{1}{\Gamma(1-\beta)}\frac{1}{t^{\beta-\alpha}}
(\www{A}^{-1}\www{a_N})(x_0) = 0                \eqno{(3.13)}
$$
by $\beta > \alpha$.

Therefore, applying (3.11) - (3.13) in (3.10) and letting 
$t \to \infty$, we obtain
$$
\left\vert \frac{1}{\Gamma(1-\alpha)}(A^{-1}a_N)(x_0) \right\vert
\le C\rho_N.                \eqno{(3.14)}
$$
Since $a_N \longrightarrow a$ in $L^2(\OOO)$ as $N\to \infty$, we see that 
$A^{-1}a_N \longrightarrow A^{-1}a$ in $\DDD(A)$, and the Sobolev
embedding implies 
$$
\lim_{N\to\infty} \Vert A^{-1}a_N - A^{-1}a\Vert_{C(\ooo{\OOO})}
= 0,
$$
that is, 
$$
\lim_{N\to\infty} A^{-1}a_N(x_0) = A^{-1}a(x_0).
$$
Hence, letting $N\to \infty$ in (3.14), we see 
$$
A^{-1}a(x_0) =  0.
$$
We set $f:= A^{-1}a$ in $\OOO$.  Then $Af = a$ in $\OOO$.
By (1.10), without loss of generality, we can assume that 
$a \le 0$ on $\ooo{\OOO}$.  In view of $c\le 0$ from (1.8), by noting
that $f\vert_{\ppp\OOO} = 0$, the weak maximum principle 
(e.g., Theorem 3.1 (p.32) in Gilbarg and Trudinger \cite{GT}) yields
$f(x) \le 0$ for $x\in \ooo{\OOO}$.

Since $f\le 0$ on $\ooo{\OOO}$ and $f(x_0) = A^{-1}a(x_0) = 0$, 
we see that $f$ attains the maximum $0$
at an interior point $x_0\in \OOO$.  In view of $(-A)f \ge 0$ in $\OOO$,
we can apply the strong maximum principle (e.g., Theorem 3.5 (p.35) in 
\cite{GT}) to conclude that $f(x)$ is a constant function.
Since $f\in \DDD(A)$, we see $f\vert_{\ppp\OOO} = 0$, so that 
$f=0$ in $\OOO$.  Then $u=0$ in $\OOO\times (0,\infty)$.
In terms of (3.8) and $\www{u}(x_0,t) = u(x_0,t) = 0$ for $t>0$, we have
$$
\frac{1}{\Gamma(1-\beta)t^{\beta}} (\www{A}^{-1}\www{a_N})
(x_0) + \www{S_N}(t) + \www{W_N}(t) = 0, \quad t\ge t_0.
$$
Multiplying with $t^{\beta}$ and letting $t\to\infty$, we similarly 
obtain 
$$
\left\vert \frac{1}{\Gamma(1-\beta)}(\www{A}^{-1}\www{a_N})(x_0)
\right\vert \le C\www{\rho_N} \quad \mbox{for all $N\in \N$.}
$$
Letting $N \to \infty$, we see $\www{A}^{-1}\www{a}(x_0) = 0$.  
By the same way as for $a$, the weak and the strong 
maximum principles similarly yield $\www{a}=0$ in $\OOO$.
This means that $a=0$ in $\OOO$ and $\www{a} = 0$ in $\www{\OOO}$,
which is a contradiction for (1.10).
Therefore, $\alpha<\beta$ is impossible.
Similarly we can prove that $\alpha>\beta$ is impossible.
Thus the proof of Theorem is complete.
$\blacksquare$
\section{Appendix}

{\bf 4.1. Proof of (1.9).}

It sufficient to prove only for $A$.  Let $\la \in \C$ be an eigenvalue
of $A$: $A\va = \la\va$ and $\va\not\equiv 0$ in $\OOO$.
That is,
$$
\left\{ \begin{array}{rl}
&-\sumij \ppp_i(a_{ij}(x)\ppp_j\va(x)) - \sumj b_j(x)\ppp_j\va
- c(x)\va = \la \va \quad \mbox{in $\OOO$}, \\
& \va\vert_{\ppp\OOO} = 0.
\end{array}\right.
                                       \eqno{(4.1)}
$$
Henceforth let $\ooo{\la}$ denote the complex conjugate of 
$\la \in \C$.  Since $a_{ij}$, $b_j$ and $c$ are real-valued, we have
$$
-\sumij \ppp_i(a_{ij}(x)\ppp_j\ooo\va(x)) - \sumj b_j(x)\ppp_j\ooo\va
- c(x)\ooo\va = \ooo\la \ooo\va \quad \mbox{in $\OOO$}.
                                                 \eqno{(4.2)}
$$
Multiplying (4.1) and (4.2) with $\ooo\va$ and $\va$ respectively, 
integrating by parts in $\OOO$ and adding, we obtain
\begin{align*}
& \int_{\OOO} \sumij a_{ij}((\ppp_i\va)\ooo{\ppp_j\va}
+ (\ooo{\ppp_i\va})\ppp_j\va) dx 
- \int_{\OOO} b_j((\ppp_j\va)\ooo{\va} + (\ooo{\ppp_j\va})\va) dx\\
-& \int_{\OOO} 2c\vert \va\vert^2 dx 
= (\la+\ooo{\la})\int_{\OOO} \vert \va\vert^2 dx.
\end{align*}
On the other hand,
$$
((\ppp_i\va)\ooo{\ppp_j\va} + \ooo{\ppp_i\va}\ppp_j\va)
= 2(\mbox{Re}\,(\ppp_i\va)\mbox{Re}\,(\ppp_j\va)
+ \mbox{Im}\,(\ppp_i\va)\mbox{Im}\,(\ppp_j\va))
$$
and so (1.2) yields 
\begin{align*}
& \sumij a_{ij}( (\ppp_i\va)\ooo{\ppp_j\va}
+ \ooo{\ppp_i\va}\ppp_j\va) \\
=& 2\left( \sumij a_{ij}\mbox{Re}\,(\ppp_i\va)\mbox{Re}\,(\ppp_j\va)
+ \sumij a_{ij}\mbox{Im}\,(\ppp_i\va)\mbox{Im}\,(\ppp_j\va)\right)\\
\ge& 2\sigma(a_{ij}) \sumij (\vert \mbox{Re}\,(\ppp_i\va)\vert^2
+ \vert \mbox{Im}\,(\ppp_i\va)\vert^2)
= 2\sigma \vert \nabla \va\vert^2.
\end{align*}
Moreover,
\begin{align*}
& -\int_{\OOO} \sumj b_j((\ppp_j\va)\ooo{\va} + \ooo{(\ppp_j\va)}\va) dx\\
=& -\int_{\OOO} \sumj b_j\ppp_j(\vert \va\vert^2) dx 
= \int_{\OOO} \left( \sumj \ppp_jb_j\right) \vert \va\vert^2 dx.
\end{align*}
Hence, using also (1.4), we have 
\begin{align*}
& 2\sigma C(\OOO)\int_{\OOO} \vert \va\vert^2 dx
+ \int_{\OOO} \left( \left( \sumj \ppp_jb_j\right)
- 2c(x)\right) \vert \va\vert^2 dx\\
\le& 2\sigma \int_{\OOO} \vert \nabla \va\vert^2 dx 
+ \int_{\OOO} \left( \sumj \ppp_jb_j\right) \vert \va\vert^2 dx
- 2\int_{\OOO} c(x)\vert\va(x)\vert^2 dx\\
\le & (\la+\ooo{\la}) \int_{\OOO} \vert \va\vert^2 dx
= 2(\mbox{Re}\, \la) \int_{\OOO} \vert \va\vert^2 dx.
\end{align*}
In view of (1.8), we see
$$
2\sigma C(\OOO) + \sumj \ppp_jb_j(x) - 2c(x) > 0 \quad
\mbox{for all $x\in \ooo{\OOO}$}.
$$
Thus Re $\la>0$ follows and we complete the proof of (1.9).
$\blacksquare$

{\bf 4.2. Proof of Proposition 1.}

For the case of symmetric $A$, in terms of the eigenfunction expansions 
by the Mittag-Leffler functions, we can directly prove the proposition
(e.g., Sakamoto and Yamamoto \cite{SY}).  However, for non-symmetric $A$,
we have to take other way.
We make use of the argument by Vergara and Zacher \cite{VZ} which is the 
first work proving Proposition 1 (i) for a case 
where the eigenfunction expansion does not work, that is,
$A$ is symmetric but the coefficients are time dependent.
For the proof of Proposition 1, relying on Chapter 5 in \cite{KRY} which 
modifies \cite{VZ}, we here give a sketch of the proof of Proposition 1.
We can mostly repeat the argument in \cite{KRY} and so we 
describe the main differences.

For the proof, we have to return to the construction of the solution $u$ to 
(1.3) by the Galerkin method.  More precisely, let $\{\psi_k\}_{k\in \N}$ be 
an orthonormal basis in $L^2(\OOO)$ which consists of all the eigenfunctions 
of $-\Delta$ with the zero Dirichlet boundary condition.  We construct
approximating solutions to $u$ in the form of
$$
u_N(x,t):= \sum_{k=1}^N p_k^N(t)\psi_k(x), \quad N\in \N   \eqno{(4.3)}
$$
such that 
$$
u_N(x,0) = \sum_{k=1}^N (u, \psi_k)_{L^2(\OOO)}\psi_k(x)
                                                           \eqno{(4.4)}
$$
and
$$
(\pppa u_N + Au_N,\, \psi_{\ell})_{L^2(\OOO)} = 0, \quad 
1\le \ell \le N, \quad t>0,
$$
which is rewritten as
$$
\int_{\OOO} \pppa u_N(x,t)\psi_{\ell}(x) dx
+ \int_{\OOO} \sumij a_{ij}(x)(\ppp_ju_N)(x,t)\ppp_i\psi_{\ell}(x) dx
$$
$$
- \int_{\OOO} \sumj b_j(\ppp_ju_N)\psi_{\ell}(x) dx 
- \int_{\OOO} c(x)u_N(x,t)\psi_{\ell}(x) dx = 0.
                                                      \eqno{(4.5)}
$$
Equation (4.4) is equivalent to a system of time-fractional ordinary 
differential equations in $\{p_k^N(t)\}_{1\le k \le N}$.
Such $p_k^N$ exist uniquely and $u_N$ converges to 
the solution $u$ weakly in $L^2(0,T;H^2(\OOO))$.  We can refer to e.g., 
Chapter 4 in \cite{KRY}, and we omit the details. 

Multiplying (4.5) with $p_{\ell}^N(t)$ and adding over $\ell=1, ..., N$, 
and integrating by parts the third term on the 
left-hand side, we obtain
\begin{align*}
& \int_{\OOO} (\pppa u_N(x,t))u_N(x,t) dx
+ \int_{\OOO} \sumij a_{ij}(x)(\ppp_ju_N)(x,t)\ppp_iu_N(x,t) dx\\
+& \hhalf \int_{\OOO} \left( \sumj \ppp_jb_j\right) \vert u_N\vert^2 dx 
- \int_{\OOO} c(x)\vert u_N(x,t)\vert^2 dx = 0, \quad t>0.
\end{align*}
By (1.2) and (1.4), we have
\begin{align*}
& \int_{\OOO} (\pppa u_N(x,t))u_N(x,t) dx
+ C(\OOO)\sigma(a_{ij})\int_{\OOO} \vert u_N(x,t)\vert^2 dx\\
+& \hhalf \int_{\OOO} \left( \sumj \ppp_jb_j\right) \vert u_N\vert^2 dx 
- \int_{\OOO} c(x)\vert u_N(x,t)\vert^2 dx \le 0, \quad t>0.
\end{align*}
Applying (1.8), we can choose a constant $\mu_0>0$ such that
$$
\int_{\OOO} (\pppa u_N(x,t))u_N(x,t) dx
+ \mu_0\Vert u_N(\cdot,t)\Vert^2_{L^2(\OOO)} \le 0, \quad t>0
\quad \mbox{for all $N\in \N$.}
$$
Hence Lemma 5.1 in \cite{KRY} yields
$$
\Vert u_N(\cdot,t)\Vert_{L^2(\OOO)} \pppa\Vert u_N(\cdot,t)\Vert
_{L^2(\OOO)}
\le -\mu_0\Vert u_N(\cdot,t)\Vert^2_{L^2(\OOO)} \le 0, \quad t>0.
                                                \eqno{(4.6)}
$$
Then, we can formally argue as follows.  By setting $d_N(t)
:= \Vert u_N(\cdot,t)\Vert_{L^2(\OOO)}$ and assuming 
that $d_N(t) \ne 0$ for all $t>0$, inequality (4.6) yields
$$
\left\{ \begin{array}{rl}
&\pppa d_N(t) \le -\mu_0d_N(t), \quad t>0,\\
& d_N(0) = \left\Vert \sum_{k=1}^N (a, \psi_k)_{L^2(\OOO)}\psi_k
\right\Vert_{L^2(\OOO)}.
\end{array}\right.
                                                 \eqno{(4.7)}
$$
Noting that 
$$
\pppa E_{\alpha,1}(-\mu_0t^{\alpha}) = -\mu_0E_{\alpha,1}(-\mu_0t^{\alpha}),
\quad  E_{\alpha,1}(0) = 1,
$$
we apply the comparison principle (e.g., Luchko and Yamamoto \cite{LY}) to 
obtain
$$
d_N(t) \le d_N(0)E_{\alpha,1}(-\mu_0t^{\alpha}), \quad t>0.
$$
Therefore, for all $N\in \N$, we have
$$
d_N(t) \le \Vert a\Vert_{L^2(\OOO)}E_{\alpha,1}(-\mu_0t^{\alpha})
\le \frac{C\Vert a\Vert_{L^(\OOO)}}{t^{\alpha}}, \quad t>0,
$$
where we use 
$$
\vert E_{\alpha,1}(-\mu_0t^{\alpha})\vert 
\le \frac{C}{1+\mu_0t^{\alpha}}, \quad t>0
$$
(e.g., Theorem 1.6 (p.35) in \cite{Po}).  This is an outline and for example,
but we need to justify in deriving (4.7) by dividing with 
$\Vert u_N(\cdot,t)\Vert_{L^2(\OOO)}$, because it may vanish at some $t$.
Such details are still the same as the proof of Theorem 5.1 of Chapter 5
in \cite{KRY}, and so we omit.

Finally we have to prove part (ii).  Let $d=1,2,3$.  By $a \in \DDD(A)$, 
we consider
$$
\left\{ \begin{array}{rl}
& \pppa U = -AU(x,t), \quad x\in \OOO,\, t>0,\\
& U\vert_{\ppp\OOO \times (0,\infty)} = 0, \\
& U(x,0) = -Aa(x), \quad x\in \OOO.
\end{array}\right.
$$
There exists a unique solution $U \in C([0,\infty);L^2(\OOO))$ and 
applying part (i) to $U$, we have
$$
\Vert U(\cdot,t)\Vert_{L^2(\OOO)} \le \frac{C}{t^{\alpha}}
\Vert Aa\Vert_{L^2(\OOO)}, \quad t>0.
$$
Moreover the uniqueness of the solution yields $U = \pppa u$ in 
$\OOO\times (0,\infty)$.  Therefore $U = -Au$ in $\OOO\times (0,\infty)$,
and so
$$
\Vert Au(\cdot,t)\Vert_{L^2(\OOO)} \le \frac{C}{t^{\alpha}}
\Vert Aa\Vert_{L^2(\OOO)}, \quad t>0.
$$
Since $d=1,2,3$, the Sobolev embedding implies 
$$
\Vert u(\cdot,t)\Vert_{C(\ooo{\OOO})} \le \frac{C}{t^{\alpha}}
\Vert Aa\Vert_{L^2(\OOO)}, \quad t>0.
$$
Thus the proof of Proposition 1 (ii) is complete.
$\blacksquare$

\section*{Acknowledgment}
The author was supported by Grant-in-Aid for Scientific Research (S)
15H05740 and Grant-in-Aid (A) 20H00117 of 
Japan Society for the Promotion of Science and
by The National Natural Science Foundation of China
(no. 11771270, 91730303).
This paper has been supported by the RUDN University 
Strategic Academic Leadership Program.


\end{document}